\documentclass[12pt,a4paper]{article}
\usepackage{amsmath,amssymb,amsthm}

\newcommand{\kap}{\noindent \bf}

\begin{document}

\begin{center}
{\large \bf Effective Support Size} \vskip 5mm
{M.  Grend\'ar\\
KM FPV Bel University, Slovakia
\\ \it umergren@savba.sk}
\end{center}

\vspace{2mm} {\kap 1.  Effective size of support}

\vspace{2mm} 


Let $X$ be a discrete random variable which can take on values from
a finite set $\mathcal X$ of $m$ elements, with probabilities
specified by the probability mass function (pmf) $p$. The support of
$X$ is a set $\mathcal S(p(X)) \triangleq \{p: p_i > 0, i = 1, 2,
\dots, m\}$. Let $|\mathcal S(p(X))|$ denote the size of the
support.

While pmf $p = [0.5, 0.5]$ makes both outcomes equally likely, the
following pmf $q = [0.999, 0.001]$ characterizes a random variable
that can take on almost exclusively only one of two values. However,
both
$p$ and $q$ have the same size of support. 
This motivates a need for a quantity that could measure size of
support of the random variable in a different way, so that the
random variable can be according to its pmf placed in the range $[1,
m]$. We will call the new quantity/measure the effective support
size (Ess), and denote it by
$\mathbb S(p(X))$, or $\mathbb S(p)$, for short. 
The example makes it obvious that $\mathbb S(\cdot)$ should be such
that $\mathbb S(q)$ will be close to $1$, while to $p$ it should
assign value $\mathbb S(p) = 2$.


\vspace{2mm} {\kap 2. Properties of Ess}

\vspace{2mm}

Ess should have certain properties, dictated by common sense.

\smallskip

P1) $\mathbb S(p)$ should be continuous, symmetric function.

P2) $\mathbb S(\delta_m) = 1 \le \mathbb S(p_m) \le \mathbb S(u_m) =
m$; where $u_m$ denotes uniform pmf on $m$-element support,
$\delta_m$ denotes an $m$-element pmf with probability concentrated
at one point, $p_m$ denotes a pmf with $|\mathcal S(p)| = m$.

P3) $\mathbb S([p_m, 0]) = \mathbb S(p_m)$.

P4) $\mathbb S(p(X,Y)) = \mathbb S(p(X))\mathbb S(p(Y))$, if $X$ and
$Y$ are independent random variables.

\smallskip

The first two properties are obvious. The third one states that
extending support by an impossible outcome should leave Ess
unchanged. Only the fourth property needs, perhaps,  some little
discussion. Or, better, an example. Let $p(X) = [1, 1, 1]/3$ and
$p(Y) = [1, 1]/2$ and let $X$ be independent of $Y$. Then $p(X,Y) =
[1, 1, 1, 1, 1, 1]/6$.  According to P2), $\mathbb S(p(X)) = 3$,
$\mathbb S(p(Y)) = 2$ and $\mathbb S(p(X,Y)) = 6 = \mathbb
S(p(X))\mathbb S(p(Y))$. It is reasonable to require the product
relationship to hold for independent random variables with arbitrary
distributions.

The properties P1)-P4) are satisfied by $\mathbb S(p, \alpha) =
\left(\sum_{i=1}^m p_i^\alpha\right)^{\frac{1}{1 - \alpha}}$, where
$\alpha$ is positive real number, different than $1$. Note that
$\mathbb S(\cdot)$ of this form is $\exp$ of R\'enyi's entropy. For
$\alpha \rightarrow 1$, $\mathbb S(p, \alpha)$ also satisfies
P1)-P4) and takes the form of $\exp(H(p))$, where $H(p) \triangleq
-\sum_{i=1}^m p_i \log p_i$ is Shannon's entropy. 
It is thus reasonable to define $\mathbb S(p, \alpha)$ for $\alpha =
1$ this way (with the convention $0\log 0 = 0$), so that $\mathbb
S(\cdot)$ then becomes a continuous function of $\alpha$.

\vspace{2mm} {\kap 3. Selecting $\alpha$}

\vspace{2mm}

The requirements define entire class of measures of effective
support size. This opens a problem of selecting $\alpha$. 
In Table 1, $\mathbb S(p, \alpha)$ is given for various two-element
pmf's, and $\alpha = 0.001, 0.1, 0.5, 0.9, 1.0, 1.5, 2.0, 10,
\infty$. The value $\mathbb S(p, \alpha\rightarrow\infty)$  can be
found analytically.

\begin{table}
\caption{$\mathbb S(p, \alpha)$ for $\alpha = 0.001$, $0.1$, $0.5$,
$0.9$, $1.0$, $1.5$, $2.0$, $10$, $\infty$ and different $p$'s.}
\center{
\begin{tabular}{|c|c|c|c|c|c|c|}
  \hline
   &  &  & $\mathbb S(p, \alpha)$ &  &  &  \\
  $\alpha$ & $[0.5, 0.5]$ & $[0.6, 0.4]$ & $[0.7, 0.3]$ & $[0.8, 0.2]$ & $[0.9, 0.1]$ & $[1.0, 0.0]$ \\
  \hline
  0.001 & 2.000000 & 1.999959 & 1.999826 & 1.999554 & 1.998979 & 1.000000 \\
  0.1 & 2.000000 & 1.995925 & 1.982696 & 1.956233 & 1.902332 & 1.000000 \\
  0.5 & 2.000000 & 1.979796 & 1.916515 & 1.800000 & 1.600000 & 1.000000 \\
  0.9 & 2.000000 & 1.964013 & 1.856116 & 1.675654 & 1.416403 & 1.000000 \\
  1.0 & 2.000000 & 1.960132 & 1.842023 & 1.649385 & 1.384145 & 1.000000 \\
  1.5 & 2.000000 & 1.941178 & 1.777878 & 1.543210 & 1.275510 & 1.000000 \\
  2.0 & 2.000000 & 1.923077 & 1.724138 & 1.470588 & 1.219512 & 1.000000 \\
  10.0 & 2.000000 & 1.760634 & 1.486289 & 1.281379 & 1.124195 & 1.000000 \\
  $\infty$ & 2.000000 & 1.666666 & 1.428571 & 1.250000 & 1.111111 & 1.000000 \\
  \hline
\end{tabular}
}
\end{table}

From the table it can be seen that the smaller the $\alpha$, the
more $\mathbb S(\cdot, \alpha)$ ignores the actual difference
between probabilities. For $p = [0.9, 0.1]$ the difference is $0.8$,
yet $\mathbb S(p, 0.001) = 1.998979$, i.e., it interprets the pmf as
being very close to $[0.5, 0.5]$.

Based on the table, we would opt for $\mathbb S(\cdot, \alpha
\rightarrow\infty)$ as the good measure of Ess. However, for larger
$|\mathcal S|$ this choice becomes less attractive. This can be seen
easily from a consideration of continuous random variables.


\vspace{2mm} {\kap 4. Selecting $\alpha$: continuous case}

\vspace{2mm}

In the case of continuous random variable $\mathbb S(f(x), \alpha)
\triangleq (\int f^\alpha(x) dx)^{\frac{1}{1-\alpha}}$. For gaussian
$n(\mu, \sigma^2)$ distribution, $\mathbb S(\cdot, \alpha) =
\frac{\sqrt{2\pi\sigma^2}}{\alpha^{\frac{1}{2(1 - \alpha)}}}$; cf.
[3]. This for $\alpha \rightarrow \infty$ converges to
$\sqrt{2\pi\sigma^2}$, so that for $\sigma^2 = 1$ it becomes
$\sqrt{2\pi} = 2.5067$. It is worth comparing with $\mathbb
S(\cdot,\alpha = 1) = \sqrt{2e\pi\sigma^2}$ (cf. [1]), which reduces
in the case of $\sigma^2 = 1$ to $4.1327$. This makes much more
sense.

That $\mathbb S(\cdot, \alpha \rightarrow \infty)$ is not the
appropriate measure of Ess can be even more clearly seen in the case
of the Exponential distribution. For $\beta e^{-\beta x}$ with
$\beta = 1$, $\mathbb S(\cdot, \alpha \rightarrow \infty) = 1$ while
$S(\cdot,\alpha = 1) = e$.

\vspace{2mm} {\kap 5. Adding another property}

\vspace{2mm}

The above considerations suggest that $\mathbb S(\cdot, \alpha = 1)$
might be the most appropriate of the Ess measures which satisfy the
requirements P1)-P4). The question is whether there is some other
requirement that is reasonable to add to the already employed
properties, such that it could narrow down the set of feasible
$\mathbb S(\cdot, \alpha)$ to $\mathbb S(\cdot, 1)$.

To this end, let us consider two random variables $X$, $Y$ that are
dependent. Let $p(Y|X)$ be the conditional distribution and $p(X,Y)$
the joint distribution. For any of them its Ess can be obtained by
$\mathbb S(\cdot, \alpha)$. For instance, let $X$ can take on two
values $x_1, x_2$. Then, $\mathbb S(p(Y|X=x_1), \alpha)$ is Ess of
the conditional distribution of $Y$ given that $X$ has taken the
value $x_1$.

In analogy with P4) it seems reasonable  to define Ess for a mean of
the conditional distributions $\mathbb S(\overline{p}(Y|X), \alpha)$
as $\mathbb S(\overline{p}(Y|X), \alpha) \triangleq \frac{\mathbb
S(p(X, Y), \alpha)}{\mathbb S(p(X), \alpha)}$. Note that $\mathbb
S(\overline{p}(Y|X), \alpha)$ is the same regardless of what value
the conditioning variable $X$ has taken. This is why it is a kind of
Ess for a mean of the conditional distributions. Note also that when
$X$ and $Y$ are independent the definition reduces to the
requirement P4).

Now, once the new object is defined, one might wonder whether it can
be related to Ess's of the conditional distributions. For $\alpha =
1$ such a relationship indeed exists:
\begin{equation}
 \mathbb S(\overline{p}(Y|X), 1) = \prod_{i=1}^n \mathbb S(p(Y|X =
 x_i), 1)^{p(X = x_i)}.
\end{equation}

If Eq. (1) was turned into the fifth requirement, then by invoking
Khinchin's [2] uniqueness theorem (which characterizes Shannon's
entropy), it can be claimed that $\mathbb S(\cdot, 1)$ is the only
Ess which satisfies the enhanced set of requirements.

It should be added, however, that Eq. (1) is not the only
perceivable relationship between $S(\overline{p}(Y|X)$ and Ess's for
conditional distributions. Instead of the form of weighted geometric
mean the relationship could for instance take the form of weighted
arithmetic mean. Whether in this case there is some $\alpha$ which
could satisfy the relationship remains to be an open problem (at
least for the present author).

\vspace{2mm} {\kap 6. Summary}

\vspace{2mm}

In this speculation we entertained the newly-introduced concept of
effective support size (Ess). There are some obvious requirements
P1)-P4) that Ess has to satisfy. The class of Ess measures $\mathbb
S(\cdot. \alpha) = \left(\sum_{i=1}^m p_i^\alpha\right)^{\frac{1}{1
- \alpha}}$ which satisfies the requirements is broad. The Ess
measures are in a direct relationship to the family of R\'enyi's
entropies which includes as its special case also Shannon's entropy.
We have briefly addressed the issue of selecting $\alpha$ such that
the corresponding $\mathbb S(\cdot, \alpha)$ would be  the most
'appropriate' measure of Ess. The considerations indicate that
$\alpha = 1$ could, perhaps, be the most reasonable candidate. If
Eq. (1) was added into the set of requirements, then $\mathbb
S(\cdot,1)$ would become the only $\mathbb S(\cdot)$ that satisfies
them. However, there are also other conceivable relationships
between $S(\overline{p}(Y|X)$ and the conditional $\mathbb
S(\cdot)$. Whether some of them could be satisfied by $\mathbb
S(\cdot, \alpha)$ for some other $\alpha$ remains to be an open
question. In any case, with the concept of Ess it is possible to
enter a meaningful world which is in a sense dual to that of
entropies.

\vspace{2mm} {\kap References}

\begin{enumerate}
 \item Verdugo Lazo, A. C. G., Rathie, P. N.: On the entropy of
 continuous probability distributions, IEEE Trans. IT, IT:24 (1978),
 pp. 120-122.
 \item Khinchin, A. I.: Mathematical foundations of Information
 Theory, Dover, NY, 1957.
 \item Song, K.-S.: R\'enyi information, loglikelihood and an intrinsic
 distribution measure, Jour. Stat. Inference and Planning 93 (2001),
 pp. 51-69.
\end{enumerate}

\medskip

{\bf Acknowledgement} Supported by VEGA 1/3016/06 grant. I am
grateful to Michael George for posing a problem which has induced
this work.

\bigskip
\medskip

\rightline{\emph{To George}.  May 2, 2006.}

\end{document}